\documentclass[10pt]{article}
\font\smallit=cmti10
\font\smalltt=cmtt10

\usepackage{amssymb,latexsym,amsmath,epsfig,amsthm} 

\makeatletter

\renewcommand\section{\@startsection {section}{1}{\z@}
{-30pt \@plus -1ex \@minus -.2ex}
{2.3ex \@plus.2ex}
{\normalfont\normalsize\bfseries\boldmath}}

\renewcommand\subsection{\@startsection{subsection}{2}{\z@}
{-3.25ex\@plus -1ex \@minus -.2ex}
{1.5ex \@plus .2ex}
{\normalfont\normalsize\bfseries\boldmath}}

\renewcommand{\@seccntformat}[1]{\csname the#1\endcsname. }

\makeatother

\newtheorem{conjecture}{Conjecture}

\newtheorem{remark}{Remark}
\newtheorem{observation}{Observation}

\begin{document}

\begin{center}
\uppercase{\bf Hecke groups, linear recurrences, and Kepler limits}
\vskip 20pt
{\bf Barry Brent}\\
{\tt barrybrent@iphouse.com}\\ 
\vskip 10pt
\end{center}
\vskip 20pt
\centerline{\smallit Received: , Revised: , Accepted: , Published: } 
\vskip 30pt

\centerline{\bf Abstract}

\noindent
We study the linear fractional
transformations  in
the Hecke group $G(\Phi)$
where $\Phi$ is either root
of  $x^2 - x -1$
(the larger root being
the ``golden ratio''
$\phi = 2 \cos \frac {\pi}5$.)
Let $g \in G(\Phi)$ and let $z$ be
a generic element of the upper half-plane.
Exploiting the fact that
$\Phi^2 = \Phi + 1$,
we find that
$g(z)$ is a quotient of linear
polynomials in $z$
such that
the coefficients
of $z^1$ and $z^0$
in the numerator and denominator of
$g(z)$
appear themselves to be linear polynomials
in $\Phi$ with coefficients
that are certain multiples of Fibonacci numbers.
 We make somewhat less
 detailed observations
 along similar lines
 about  the functions
 in $G(2 \cos \frac {\pi}k)$
 for $k$ greater than or equal to $5$.

\pagestyle{myheadings} 
\markright{\smalltt INTEGERS: 19 (2019)\hfill} 
\thispagestyle{empty} 
\baselineskip=12.875pt 
\vskip 30pt
\section{Introduction}
Let $G(\lambda)$ be the
Hecke group
 generated by
the linear fractional transformations
$S \colon z \mapsto -1/z$ and
$T_{\lambda} \colon z \mapsto z + \lambda$
and let $G_k = G(2 \cos \frac {\pi}k)$.
This article describes numerical experiments
carried out
to study Hecke groups,
mainly
$G_k$ for $k \geq 5$.
In this article, an $n$-tuple of symbols
$$\overrightarrow{w} = \{w_1, w_2, w_3, ..., w_n\}$$
representing an ordered set of integers
is called a \it word \rm on $\bf{Z}$
and we write $|\overrightarrow{w}|$
for $n$.
A typical
element of $G(\lambda)$
takes the form
$$T_{\lambda}^{w_1}ST_{\lambda}^{w_2}S ... S
T_{\lambda}^{w_n} =
 g_{_{\small{\lambda}};
_{\overrightarrow{w}}}.$$

This representation is not unique.
For example,
 a function
$g$ in $G(\lambda)$
can be described by a word
of length $n$ for
arbitrarily large $n$,
because
any word representing $g$
can be padded with
zeros and the resulting word will
also represent $g$. Consequently, when
studying all $g$ represented by words
$\overrightarrow{w}$
with $|\overrightarrow{w}|$ 
less than or equal to $N$,
we can restrict attention to
the words $w$ such that
$|\overrightarrow{w}|$ is
equal to $N$.

Let $\phi, \phi^* $,
represent the larger and smaller roots 
of $x^2 - x - 1$,
respectively.
The problem of expressing  
(for $g \in G_5)$
$g = g_{_{\phi};_{\small{\overrightarrow{w}}}}$
 in terms of the $w_i$
was raised by Leo in [9]
and discussed by his student Sherkat in
[12];
the first purpose of this article is to
write down conjectures addressing this question.
Our calculations indicate that,
for arbitrary $\lambda, z \in \bf{C}$,
$g_{_{\small{\lambda}};_{\overrightarrow{w}}}(z)$
is a rational function of $z$ and $\lambda$  in
polynomials of
$\lambda$-degree
less than or equal to
$k$. Here are the first few:

$$g_{\small{\lambda};\{w_1\}}(z) =
\frac{1\cdot z + w_1 \lambda}{0 \cdot z + 1},
$$

$$
g_{\small{\lambda};\{w_1,w_2\}}(z) =
\frac{w_1 \lambda \cdot z + w_1 w_2 \lambda^2 - 1}
{1 \cdot z +w_2 \lambda},
$$

and
$$
g_{\small{\lambda};\{w_1,w_2,w_3\}}(z) =
\frac{(w_1 w_2\lambda^2 - 1) \cdot z +
w_1 w_2 w_3 \lambda^3 - (w_1 + w_3)\lambda}
{w_2 \lambda \cdot z + w_2 w_3 \lambda^2 -1}.
$$
Following  [12],
we simplify the above  expressions for
$g_{_{\small{\lambda}};_{\small{\overrightarrow{w}}}}$
when
$\lambda$ is $\phi$ or $\phi^*$
by repeatedly making the substitution
$\Phi^2 = \Phi + 1$ ($\Phi = \phi$ or 
$\Phi = \phi^*$.)
The coefficients in
$g_{_{\small{\Phi};_{\small{\overrightarrow{w}}}}}$
become linear polynomials
in $\Phi$:
\newline \newline

$$g_{\small{\Phi};\{w_1\}}(z) =\frac{1 \cdot z
+ w_1 \Phi}{0 \cdot z + 1},
$$

$$g_{\small{\Phi};\{w_1,w_2\}}(z) =
\frac{w_1 \Phi \cdot z + w_1 w_2 \Phi + w_1 w_2 - 1}
{1\cdot z +w_2 \Phi},
$$

and
$$g_{\small{\Phi};\{w_1,w_2,w_3\}}(z) = $$
$$\frac{(w_1 w_2\Phi + w_1 w_2 - 1) \cdot z
+ (2w_1 w_2 w_3 -w_1 - w_3) \Phi +
w_1 w_2 w_3}
{w_2 \Phi \cdot z + w_2 w_3 \Phi + w_2 w_3 -1}.
$$
Further calculations suggest that
the coefficients of $\Phi^1$ and $\Phi^0$ in
these expressions are always
linear combinations of first-degree monomials
$h$
in the $w_i$ such that
the numerical coefficient
of $h$ is $\pm 1$ times a Fibonacci number
determined by the total degree of $h$;
details are in the next section.

It is well known, of course, that the
consecutive ratios $F_n/F_{n-1}$
of Fibonacci numbers
converge to $\phi$.
More generally,
the limit of the ratios of
consecutive elements of a linear recurrence $L$,
when it exists, is called by Fiorenza and Vincenzi
the \it Kepler limit \rm of $L$.
Certain roots of polynomials other
than $x^2 - x -1$ are also Kepler limits
[5, 6], so we are led to consider
the possibility that the $G_5$ phenomenon generalizes;
our observations tend to confirm this guess.
Section 2 describes what we found
out about $G_5$,
Section 3 describes less detailed observations for
$G_k$ with $k$ between $5$ and $33$ (inclusive),
and the final section provides some detail about
our numerical experiments;
documentation in the form of \it Mathematica \rm
notebooks is on our ResearchGate site for 
this article [4].

We state merely
 empirical claims in
this article.
We make several conjectures,
but they, too, are based 
on empirical evidence,
not on 
theoretical reasoning. 
When we say we have observed convergence
of a sequence $s_n$ (say)
to a limit $S$, we mean that
our plots of $1000$
values of $\log |S-s_n|$ are
apparently linear, with negative slope.
We rely on our eyesight in this matter: 
we have not fitted our data to curves with
a statistical package. Interested readers are
invited to inspect the plots
on our ResearchGate pages.

In the following section our observations
were made on words
in $W$ of length $20$,
and those in the last section
tested words of length $25$.
This means that we have in fact
tested the claims
on all shorter words as well.

In our tests, we identified
functions in the $G_k$ with 
their matrix representations:
a function
$$T_{\lambda}^{w_1}ST_{\lambda}^{w_2}S ... S
T_{\lambda}^{w_n}$$
was identified
with the corresponding matrix product.

More information about the
Hecke groups is available,
for example,
in [2].
\begin{remark} \rm
The book [7]
by Khovanskii
apparently describes another method for
approximating roots
of polynomials
using convergent sequences
of ratios of elements
of numerical sequences; but these
sequences are not linear recurrences.
(We have
not seen [7],
but Khovanskii's's method
is described in [10], where
the book is cited.)
\end{remark}
\section[]{The group {$G_5$}}
We make the following definitions.
\newline\newline\noindent
1. The Fibonacci numbers are defined with the convention that
$F_0 = 0, F_1 = F_2 = 1$, \it etc. \rm It will be convenient
to write $F_{-1} = 1$ as well in contexts where (see below)
$\overrightarrow{s} = \emptyset$.
\newline\newline\noindent
2. $\chi$ is the following Dirichlet character:
\[
\chi(n) =
\left\{
\begin{array}{cl}
0 & \mbox{if $n \equiv 0 \pmod{4}$,}\\
1 & \mbox{if $n \equiv 1 \pmod{4}$,}\\
0 & \mbox{if $n \equiv 2 \pmod{4}$,}\\
-1 & \mbox{if $n \equiv 3 \pmod{4}.$}\\
\end{array}
\right.
\]
Alternatively, with $(a|b)$ representing
the Kronecker symbol,
$\chi(n) = (n|2)$ if
$n \equiv 0, 1, 2, 3, 4$ or $6$ $\pmod{8}$,
and $\chi(n) = -(n|2)$ otherwise.
\newline\newline\noindent
3. $W$ is the set of words $\overrightarrow{w}$ on $\bf{Z}$.
The empty word $\overrightarrow{\emptyset}$
verifies
$\overrightarrow{\emptyset} \in W$
and $\overrightarrow{\emptyset} \subset \overrightarrow{w}$
for any $\overrightarrow{w} \in W$.
\newline\newline\noindent
4. We write the
cardinal number of a set $\sigma$
as $|\sigma|$.
We apply the same notation to words 
$\overrightarrow{w}$ in $W$.
We write $|\overrightarrow{\emptyset}| = 0$.
\newline\newline\noindent
5. (a) For 
$\overrightarrow{w} \in W, \overrightarrow{w} 
\neq \overrightarrow{\emptyset}$,
let
$\overrightarrow{s} = \{w_{j_1}, w_{j_2}, ..., w_{j_m}\}
\subset \overrightarrow{w} = \{w_1, ..., w_n \}$.
If all of the $j_m \equiv m \pmod{2}$,
\newline\newline\noindent
then we  write
$$\overrightarrow{s} \ll _1 \overrightarrow{w}.$$
We also write
$\overrightarrow{\emptyset} \ll _1 \overrightarrow{w}$.
\newline\newline\noindent
 (b) If
 $\overrightarrow{s} \ll _1 \overrightarrow{w}$
 and
$|\overrightarrow{s}| > 1$,
we write $$\overrightarrow{s} \ll _2 \overrightarrow{w}.$$
\noindent
(c) Let $\overrightarrow{s}, \overrightarrow{w}$ be as
in definition 5a,
except that all of the $j_m$ satisfy
$j_m \equiv m - 1 \pmod{2}$.
Then we  write $$\overrightarrow{s} \ll _3 \overrightarrow{w}.$$
We also write
$\overrightarrow{\emptyset} \ll _3 \overrightarrow{w}$.
\newline\newline\noindent
(d) If
$\overrightarrow{s} \ll _3 \overrightarrow{w}$
and $|\overrightarrow{s}| > 1$,
we write $\overrightarrow{s} \ll _4 \overrightarrow{w}$.
\newline\newline\noindent
6. (a) For
 $\overrightarrow{w} \in W$, the formal product
$$m_{\small{\overrightarrow{w}}} =
\prod_{w_i \in \small{\overrightarrow{w}}} w_i .$$
We also write
$$m_{\small{\overrightarrow{\emptyset}}} = 1.$$
\newline\newline\noindent
(b) $M_{\overrightarrow{w}}$ is the set of all
linear combinations with coefficients in the integers
of monomials $m_{\small{\overrightarrow{s}}}$
such that $\overrightarrow{s} \subset \overrightarrow{w}$.
\newline\newline\noindent
(c) $M$ is the union of the $M_{\overrightarrow{w}}$
as $\overrightarrow{w}$ ranges over $W$.
\begin{remark} \rm  
In view of the identities $\Phi^2 = \Phi + 1$
for $\Phi = \phi$ or $\phi^*$,
it is clear that
\newline\newline\noindent
(i) For each $j$ between $1$ and $8$ (inclusive), 
there is a function
$f_j \colon W \rightarrow M$ such that
$f_j(\overrightarrow{w}) \in M_{\small{\overrightarrow{w}}}$
and,
for all $g_{_{\small{\Phi};
_{\small{\overrightarrow{w}}}}} \in G(\Phi)$ 
and $z$ with
$\Im z$ positive,
$$g_{_{\small{\Phi};_{\small{\overrightarrow{w}}}}} (z) =
\frac{(f_1(\small{\overrightarrow{w}})
\Phi + f_2(\small{\overrightarrow{w}})) z
+
f_3(\small{\overrightarrow{w}})\Phi +f_4(\small{\overrightarrow{w}})
}{(f_5(\small{\overrightarrow{w}})
\Phi + f_6(\small{\overrightarrow{w}})) z
+
f_7(\small{\overrightarrow{w}})\Phi +
f_8(\small{\overrightarrow{w}})}.
$$
Referring to the introduction, for example:
$$
f_3(w_1, w_2, w_3) =
 2w_1 w_2 w_3 -w_1 - w_3
$$
and
$$f_6(w_1, w_2, w_3) = 0.$$
\noindent
(ii) For each $j$ between $1$ and $8$
(inclusive),
there is a function
$\nu_j \colon  W \times W \mapsto \bf{Z}$
determined by the condition
$$f_j(\small{\overrightarrow{w}}) =
\sum_{\overrightarrow{s} \subset \overrightarrow{w}}
\nu_j(\small{\overrightarrow{s}},\overrightarrow{w})
m_{\small{\overrightarrow{s}}}.$$
\end{remark}

The following observations
describe the
functions represented by words of length 
less than or equal to $20$
in $G_k$ for $k$ between $5$ and $50$ 
(inclusive.)
\begin{observation} \rm
(a) 
\[
\nu_1(\small{\overrightarrow{s}},\overrightarrow{w}) =
\left\{
\begin{array}{cl}
\chi(|\overrightarrow{w}| - |\overrightarrow{s}|)
F_{|\overrightarrow{s}| }  
& \mbox{if $\overrightarrow{s} \ll_1 \overrightarrow{w}$,}\\
0 & \mbox{otherwise.}
\end{array}
\right.
\]
(b)
\[
\nu_2(\small{\overrightarrow{s}},\overrightarrow{w}) =
\left\{
\begin{array}{cl}
\chi(|\overrightarrow{w}| - |\overrightarrow{s}|)
F_{|\overrightarrow{s}| -1}  
& \mbox{if $\overrightarrow{s} \ll_1 \overrightarrow{w}$,}\\
0 & \mbox{otherwise.}
\end{array}
\right.
\]
(c)
\[
\nu_3(\small{\overrightarrow{s}},\overrightarrow{w}) =
\left\{
\begin{array}{cl}
-\chi(|\overrightarrow{w}| - |\overrightarrow{s}|-1)
F_{|\overrightarrow{s}|}  
& \mbox{if $\overrightarrow{s} \ll_1 \overrightarrow{w}$,}\\
0 & \mbox{otherwise.}
\end{array}
\right.
\]
(d)
\[
\nu_4(\small{\overrightarrow{s}},\overrightarrow{w}) =
\left\{
\begin{array}{cl}
-\chi(|\overrightarrow{w}| - |\overrightarrow{s}|-1)
F_{|\overrightarrow{s}|-1 }  
& \mbox{if $\overrightarrow{s} \ll_2 \overrightarrow{w}$,}\\
0 & \mbox{otherwise.}
\end{array}
\right.
\]
(e)
\[
\nu_5(\small{\overrightarrow{s}},\overrightarrow{w}) =
\left\{
\begin{array}{cl}
\chi(|\overrightarrow{w}| - |\overrightarrow{s}| - 1)
F_{|\overrightarrow{s}|-1 }  
& \mbox{if $\overrightarrow{s} \ll_3 \overrightarrow{w}$,}\\
0 & \mbox{otherwise.}
\end{array}
\right.
\]
(f)
\[
\nu_6(\small{\overrightarrow{s}},\overrightarrow{w}) =
\left\{
\begin{array}{cl}
\chi(|\overrightarrow{w}| - |\overrightarrow{s}| - 1)
F_{|\overrightarrow{s}|-1 }  
& \mbox{if $\overrightarrow{s} \ll_4 \overrightarrow{w}$,}\\
0 & \mbox{otherwise.}
\end{array}
\right.
\]
(g)
\[
\nu_7(\small{\overrightarrow{s}},\overrightarrow{w}) =
\left\{
\begin{array}{cl}
\chi(|\overrightarrow{w}| - |\overrightarrow{s}|)
F_{|\overrightarrow{s}|}  
& \mbox{if $\overrightarrow{s} \ll_3 \overrightarrow{w}$,}\\
0 & \mbox{otherwise.}
\end{array}
\right.
\]
(h)
\[
\nu_8(\small{\overrightarrow{s}},\overrightarrow{w}) =
\left\{
\begin{array}{cl}
\chi(|\overrightarrow{w}| - |\overrightarrow{s}|)
F_{|\overrightarrow{s} | - 1}  
& \mbox{if $\overrightarrow{s} \ll_3 \overrightarrow{w}$,}\\
0 & \mbox{otherwise.}
\end{array}
\right.
\]
\end{observation}
\begin{conjecture}
Observation 1
holds for words of arbitrary length and 
all $k$ greater than or equal to $5$.
\end{conjecture}
\section{Higher-order Hecke groups}
\bf{Definition}\rm:
Let  $t(x)$ be a polynomial
$ \sum_{j=0}^d a_j x^j$
and
$\gamma(t) : = \gcd \{j $   s.t. $a_j \neq 0\} $.
If $\gamma(t) = 1$,
we say that $t$ is \it stable \rm.
Whether or not $t$
is stable, we  associate to it
the family of $d^{th}$-order linear recurrences
$\Lambda_{t}$
with kernel
$\{-a_{d-1}, -a_{d-2}, ..., -a_0 \}$.
Let $\lambda = 2 \cos \frac{\pi}k$
with
minimal polynomial
$p_{\lambda} = p$ (say.)
Under certain conditions
[5, 6],
a
root $x = \kappa_p$ of $p(x)$
is the Kepler limit of one of the
$L_p \in \Lambda_p$.
The elements
of $G(\lambda) = G_k$
have the form
\begin{equation}
g_{_{\small{\lambda};_{\small{\overrightarrow{w}}}}} (z)
= \frac{(\sum_{j=0}^{d-1} f_{\lambda,1,j}
(\small{\overrightarrow{w}})\lambda^j)\cdot z
+ \sum_{j=0}^{d-1}
f_{\lambda,2,j}(\small{\overrightarrow{w}})\lambda^j
}{(\sum_{j=0}^{d-1} f_{\lambda,3,j}
(\small{\overrightarrow{w}})\lambda^j)\cdot z
+ \sum_{j=0}^{d-1} f_{\lambda,4,j}(\small{\overrightarrow{w}})\lambda^j}.
\end{equation}
(Equation (1)
is clear, as in the $G_5$ case, by substitution.)

For pragmatic reasons, we restrict
our attention  to
$f = f_{\lambda,1,0}$ in the
following observations. 
\begin{observation}\rm
For 
$k$ between $5$ and $500$ (inclusive),
$\gamma(p) = 1$ if $k$ is odd
and
$\gamma(p) = 2$ if $k$ is even.
\end{observation}
\begin{conjecture} For polynomials
of the form $p = p_{\lambda}$,
the statements in the above observation 
hold for all $k$ greater than or equal to $5$.
\end{conjecture}
\begin{observation} \rm
Let $k$ lie between $5$ and $33$
(inclusive).
\newline \newline \noindent
 (a) There is a function
$\nu^{(k)} \colon  W \times W \mapsto \bf{Z}$ 
such that
\begin{equation}
f(\small{\overrightarrow{w}}) =
\sum_{\overrightarrow{s} \subset \overrightarrow{w}}
\nu^{(k)}(\small{\overrightarrow{s}},\overrightarrow{w})
m_{\small{\overrightarrow{s}}}.
\end{equation}
with
$$
\overrightarrow{s_1}, 
\overrightarrow{s_2} \subset \overrightarrow{w}
$$
and
\begin{equation}
|\overrightarrow{s_1}| = |\overrightarrow{s_2}| \Rightarrow
|\nu^{(k)}(\small{\overrightarrow{s_1}},\overrightarrow{w})| =
|\nu^{(k)}(\small{\overrightarrow{s_2}},\overrightarrow{w})|
\end{equation}
for all $\overrightarrow{w} \in W$
s.t. $|\overrightarrow{w}| = 25$.
\newline \newline \noindent
(b) If $k$ is odd, then for some particular
$L_p\in \Lambda_p$
and all
$\small{\overrightarrow{s}} \subset  \overrightarrow{w}$
s.t. $|\overrightarrow{w}| = 25$:
\newline \newline \noindent
(b1) $|\nu^{(k)}(\small{\overrightarrow{s}},
\overrightarrow{w})| \in L_p$ and
(b2) $\kappa_p = \lambda$.
\newline \newline \noindent
(In our experiments the sum on the right-hand side of 
Equation (2) typically contains over $6 \times 10^4$
terms, but twelve or fewer distinct values 
of 
$|\nu^{(k)}(\small{\overrightarrow{s}},\overrightarrow{w})|$.)
\newline \newline \noindent
(c) Suppose $k$ is an even number between $6$
and $32$ (inclusive.)
Then
\newline \newline \noindent
(c1) clause (b1)  still holds,
but (b2) does not; in this situation,
we found no  $L_p$ for
which  $\kappa_p$ exists.
(By design, our searches stop with the
first instance of $L_p$
satisfying (a), so this is
far from decisive.)
\newline \newline \noindent
(c2) For $k = 8, 10, 14, 16, 18, 22, 26$, and $32$,
the ratios of consecutive elements
of the $L_p$ we found in
the experiments
form two convergent sub-sequences with different limits.
\newline \newline \noindent
(c3) For $k = 6, 12, 20,  24,  28$, and $30$, the
$L_p$ terminate in a sequence in which
alternate members are zero, so that the
 requisite ratios are alternately zero or
undefined.
\newline \newline \noindent
(d)  Suppose $k =
12, 14, 20, 22, 24, 28$, or $30$.
After a substitution
$y = x^2$, $p(x)$ is
transformed to a stable polynomial
$q_{\lambda^2}(y) = q(y)$ (say),
and then
 $\lambda^2$
is the Kepler limit of
a linear recurrence 
$L_q \in \Lambda_q$
containing the $|\nu^{(k)}(\small{\overrightarrow{s}},
\overrightarrow{w})|$.
\end{observation}
\begin{conjecture}
In the above observation, clause (a)
holds for all 
$k$ greater than or equal to $5$,
clause (b) holds for all odd
$k$ greater than or equal to $5$, and
clause (c1) holds for all
even 
$k$ greater than or equal to $6$.
One of clauses (c2)  or (c3)
holds for any even $k \geq 6$.
Clause (d) holds for an unbounded
set of even 
$k$ greater than or equal to $6$.
\end{conjecture}

The conditions on polynomials
under which linear recurrences
with Kepler limits that are killed by 
them
were established in [11]
(cited in [5]).

A procedure (which can be invoked
 from computer algebra systems)
 for computing
  $p = p_{\lambda}$
 for $\lambda = 2 \cos \frac{\pi}k$
 one at a time
 for individual $k$
 appeared in
 [2];
 some information about the constant terms,
 in [1]; and,
 about the
 degree, in [8].
 \section{Data on the linear recurrences}
\subsection{Coefficients}
This is a list of distinct coefficents
of the $m_{\overrightarrow{s}}$
appearing in our calculations for 
Equation (2), 
($k$ between $5$  and $33$, inclusive):
 \newline \newline \noindent
5: 1, 2, 5, 13, 34, 89, 233, 610, 1597, 4181, 10946, 28657
 \newline \newline \noindent
6: 1, 3, 9, 27, 81, 243, 729, 2187, 6561, 19683, 59049, 177147, 531441
 \newline \newline \noindent
7: 1, 4, 14, 47, 155, 507, 1652, 5373, 17460, 56714, 184183
 \newline \newline \noindent
8: 1, 2, 8, 28, 96, 328, 1120, 3824, 13056, 44576, 152192, 519616
 \newline \newline \noindent
9: 1, 6, 27, 109, 417, 1548, 5644, 20349, 72846, 259579
 \newline \newline \noindent
10: 1, 5, 25, 100, 375, 1375, 5000, 18125, 65625, 237500, 859375, 3109375
 \newline \newline \noindent
11: 1, 6, 27, 110, 429, 1637, 6172, 23104, 86090, 319792
 \newline \newline \noindent
12: 1, 4, 15, 56, 209, 780, 2911, 10864, 40545, 151316, 564719
 \newline \newline \noindent
13: 1, 6, 27,110, 429, 1638, 6188, 23255, 87190, 326646
 \newline \newline \noindent
14: 1, 7, 49, 245, 1078, 4459, 17836, 69972, 271313, 1044435, 4002467
 \newline \newline \noindent
15: 1, 5, 20, 74, 265, 936, 3290, 11560, 40699, 143755, 509771
 \newline \newline \noindent
16: 1, 2, 16, 88, 416, 1820, 7616, 31008, 124032, 490312
 \newline \newline \noindent
17: 1, 8, 44, 208, 910, 3808, 15504, 62016, 245157
 \newline \newline \noindent
18: 1, 3, 18, 81, 333, 1323, 5184, 20196, 78489, 304722, 1182519
 \newline \newline \noindent
19: 1, 10, 65, 350, 1700, 7752, 33915, 144210
 \newline \newline \noindent
20: 1, 8, 45, 220, 1000, 4352, 18411, 76380, 312455
 \newline \newline \noindent
21: 1, 7, 35, 154, 636, 2534, 9877, 37962, 144571, 547239
 \newline \newline \noindent
22: 1, 11, 121, 847, 4840, 24684, 117249, 531069, 2326588
 \newline \newline \noindent
23: 1, 12, 90, 544, 2907, 14364, 67298
 \newline \newline \noindent
24: 1, 8, 44, 208, 911, 3824, 15656, 63136, 252241
 \newline \newline \noindent
25: 1, 10, 65, 350, 1700, 7752, 33915, 144210
 \newline \newline \noindent
26: 1, 13, 169, 1352, 8619, 48165, 247247, 1197196
 \newline \newline \noindent
27: 1, 18, 189, 1518
 \newline \newline \noindent
28: 1, 12, 91, 560, 3059, 15484, 74382
 \newline \newline \noindent
29: 1, 14, 119, 798, 4655, 24794
 \newline \newline \noindent
30: 1, 7, 35, 155, 650, 2653, 10676, 42635, 169555
 \newline \newline \noindent
31: 1, 16,152, 1120, 7084
 \newline \newline \noindent
32: 1 , 2, 32, 304, 2240, 14168
 \newline \newline \noindent
33: 1, 11, 77, 440, 2244, 10659, 48278, 211486
\subsection{Initial segments for observations 3a - 3c}
This section
describes the results of a search
for initial segments $I$ of linear recurrences 
$L_p (p = p_{\lambda}, \lambda = 2 \cos \frac {\pi}k)$
with length equal to that of the kernel of
$L_p$ (so that $I$ determines $L_p$)
such that a sufficiently long initial segment of
$L_p$ contains the elements listed above 
for corresponding $k$. 
\newline \newline \noindent
5: \{0, 1\} \newline \newline \noindent
6: \{0, 1\} \newline \newline \noindent
7: \{0, 0, 1\} \newline \newline \noindent
8: \{0, 0, 1, 2\} \newline \newline \noindent
9: \{0, 0, 1\} \newline \newline \noindent
10: \{0, 0, 1, 5\} \newline \newline \noindent
11: \{0, 0, 0, 0, 1\} \newline \newline \noindent
12: \{0, 0, 0, 1\} \newline \newline \noindent
13: \{0, 0, 0, 0, 0, 1\} \newline \newline \noindent
14: \{ -3, 1, -3, 0, -3, 0 \} \newline \newline \noindent
15: \{0, 0, 0, 1\} \newline \newline \noindent
16: \{0, 0, 0, 0, 0, 0, 1, 2\} \newline \newline \noindent
17: \{0, 0, 0 ,0, 0, 0, 0, 1\} \newline \newline \noindent
18: \{0, 0, 0, 0, 1, 3\} \newline \newline \noindent
19: \{0, 0, 0, 0, 0, 0, 0, 0, 1\} \newline \newline \noindent
21: \{0, 0, 0, 0, 0, 1\} \newline \newline \noindent
22: \{-1, 1, -1, 0, -1, 0, -1, 0, -1, 0\} \newline \newline \noindent
23: \{0, 0, 0, 0, 0, 0, 0, 0, 0, 0, 1\} \newline \newline \noindent
24: \{0, 0, 0, 0, 0, 0, 0, 1\} \newline \newline \noindent
25: \{0, 0, 0, 0, 0, 0, 0, 0, 0, 1\} \newline \newline \noindent
26: \{-1, -1, -1, 0 , -1, 0, -1, 0, -1, 0, 1, 0\} \newline \newline \noindent
27: \{0, 0, 0, 0, 0, 0, 0, 0, 1\} \newline \newline \noindent
28: \{0, 0, 0, 0, 0, 0, 0, 0, 0, 0, 0, 1\} \newline \newline \noindent
29: \{0, 0, 0, 0, 0, 0, 0, 0, 0, 0, 0, 0, 0, 1\} \newline \newline \noindent
30: \{0, 0, 0, 0, 0, 0, 0, 1\} \newline \newline \noindent
31: \{0, 0, 0, 0, 0, 0, 0, 0, 0, 0, 0, 0, 0, 0, 1\} \newline \newline \noindent
32: \{0, 0, 0, 0, 0, 0, 0, 0, 0, 0, 0, 0, 0, 0, 1, 2\} \newline \newline
\noindent
33: \{0, 0, 0, 0, 0, 0, 0, 0, 0, 1\}
\subsection{Initial segments for observation 3d}
This is a list of initial segments for
$L_q \in \Lambda_q,
q=q_{\lambda^2}, \lambda = 2 \cos \frac {\pi}k, 
k= 12, 14, 20, 22, 24, 28$, and $30$, satisfying the conditions
of observation 3d.
\newline \newline \noindent
12: \{0, 1\} \newline \newline \noindent
14: \{1, 0, 0\}  \newline \newline \noindent
20: \{0, 0, 1\}  \newline \newline \noindent
22: \{1, 0, 0, 0, 0\}  \newline \newline \noindent
24: \{0, 0, 0, 1 \}  \newline \newline \noindent
28: \{0, 0, 0, 0, 0, 1 \}  \newline \newline \noindent
30: \{0, 0, 0, 1\}  \newline \newline \noindent
\subsection{Kernels for the linear recurrences}
We list these for the convenience of the
reader. Below is the list of kernels for the 
$L_p \in \Lambda_p, (p = p_{\lambda},
\lambda = 
2\cos \frac {\pi}k$ 
for $k$ between $5$ and $33$
(inclusive.)
 \newline \newline \noindent
5: \{1 ,1\}\newline \newline \noindent
6: \{0, 3\}\newline \newline \noindent
7: \{1, 2, -1\}\newline \newline \noindent
8: \{0, 4, 0, -2\}\newline \newline \noindent
9: \{0, 3, 1\}\newline \newline \noindent
10: \{0, 5, 0, -5\}\newline \newline \noindent
11:  \{1, 4, -3, -3, 1\}\newline \newline \noindent
12:  \{0, 4, 0, -1\}\newline \newline \noindent
13:  \{1, 5, -4, -6, 3, 1\}\newline \newline \noindent
14: \{0, 7, 0, -14, 0, 7\}\newline \newline \noindent
15: \{-1, 4, 4, -1\}\newline \newline \noindent
16: \{0, 8, 0, -20, 0, 16, 0, -2\}
\newline \newline \noindent
17: \{1, 7, -6, -15, 10, 10, -4, -1\}
\newline \newline \noindent
18: \{0, 6, 0, -9, 0, 3\}\newline \newline \noindent
19: \{1, 8, -7, -21, 15, 20, -10, -5, 1\}
\newline \newline \noindent
20: \{0, 8, 0, -19, 0, 12, 0, -1\}
\newline \newline \noindent
21: \{-1, 6, 6, -8, -8, -1\}\newline \newline\noindent
22: \{0, 11, 0, -44, 0, 77, 0, -55, 0, 11\}
\newline \newline \noindent
23: \{1, 10, -9, -36, 28, 56, -35, -35, 15, 6, -1\}
\newline \newline \noindent
24: \{0, 8, 0, -20, 0, 16, 0, -1\}
\newline \newline \noindent
25: \{0, 10, 0, -35, 1, 50, -5, -25, 5, 1\}
\newline \newline \noindent
26: \{0, 13, 0, -65, 0, 156, 0, -182, 0, 91, 0, -13\}
\newline \newline \noindent
27: \{0, 9, 0, -27, 0, 30, 0, -9, 1\}
\newline \newline \noindent
28: \{0, 12, 0, -53, 0, 104, 0, -86, 0, 24, 0, -1\}
\newline \newline \noindent
29: \{1, 13, -12, 
-66, 55, 165, -120, -210, 126, 126, -56, -28, 7, 1\}
\newline \newline \noindent
30: \{0, 7, 0, -14, 0, 8, 0, -1\}
\newline \newline \noindent
31: \{1, 14, -13, -78, 66, 220, -165, 
-330, 210, 252, -126, -84, 28, 8, -1
\}\newline \newline \noindent
32: \{0, 16, 0, -104, 0, 352, 0,
-660, 0, 672, 0, -336, 0, 64, 0, -2\}
\newline \newline \noindent
33: \{-1, 10, 10, -34, -34, 43, 43, -12, -12, -1\}
\newline \newline \noindent
Below is the list of kernels for the 
$L_q \in \Lambda_q,
q=q_{\lambda^2}, \lambda = 2 \cos \frac {\pi}k, 
k= 12, 14, 20, 22, 24, 28$, and $30$, satisfying the conditions
of observation 3d.
\newline \newline \noindent
12: \{4, -1\}
\newline \newline \noindent
14: \{7, -14, 7\}
\newline \newline \noindent
20: \{8, -19, 12, -1\}
\newline \newline \noindent
22: \{11, -44, 77, -55, 11\}
\newline \newline \noindent
24: \{8, -20, 16, -1\}
\newline \newline \noindent
28: \{12, -53, 104, -86, 24, -1\}
\newline \newline \noindent
30: \{7, -14, 8, -1\}
\begin{center}
    \bf{References}
\end{center}
[1] C. Adiga, I. N. Cangul and H. N. Ramaswamy, On the constant term 
of the minimal polynomial of $\cos \frac {2\pi}n$ over $\mathbb{Q}$,
{\it Filomat}, {\bf 30}  (2016), 1097--1102.
\newline \newline \vskip -.25in \noindent
[2] A. Bayad and I. N. Cangul, The minimal polynomial of
$2\cos (\frac {\pi}q)$ and Dickson polynomials, 
{\it Appl. Math. Comp.}, {\bf 218} 
(2012), 7014--7022.
\newline \newline \vskip -.25in \noindent
[3] B. C. Berndt and M. I. Knopp,
{\it Hecke's Theory of Modular Forms and Dirichlet Series},
World Scientific, Singapore, 2008.
\newline \newline \vskip -.25in \noindent
[4] B. Brent, https://www.researchgate.net/profile/
\newline (search under ``Barry Brent'', including the quotation marks.)
\newline\newline \vskip -.25in \noindent
[5] A. Fiorenza and G. Vincenzi, From Fibonacci sequence to
the golden ratio, \newline http://dx.doi.org/10.1155/2013/204674
{\it J. Math.} {\bf 2013} (2013) 
\newline \newline \vskip -.25in \noindent
[6] A. Fiorenza and G. Vincenzi, 
Limit of ratio of consecutive terms for general order-k linear 
homogeneous recurrences with constant coefficients,
{\it Chaos Solitons Fractals} {\bf 44} (2011), 145--152.
\newline \newline \vskip -.25in \noindent
[7] A. N. Khovanskii, {\it The Application of Continued Fractions
and Their Generalizations to Problems in Approximation Theory}
(tr. P. Wynn), Noordhoff Publishers, Groningen, 1963.
\newline \newline \vskip -.25in \noindent
[8] D. H. Lehmer, A note on trigonometric algebraic numbers,
{\it Amer. Math. Monthly}, {\bf 40} (1933), 165-166.
\newline \newline \vskip -.25in \noindent
[9] J. G. Leo, {\it Fourier Coefficients of Triangle  Functions}
(Ph. D. thesis), \newline http://halfaya.org/ucla/research/thesis.pdf,
U.C.L.A, Los Angeles, 2008.
\newline \newline \vskip -.25in \noindent
[10] J. Mc Laughlin and B. Sury, Some observations on Khovanskii's
matrix methods for extracting roots in polynomials, 
{\it Integers}, {\bf 7} (2007), \# A48.
\newline \newline \vskip -.25in \noindent
[11] H. Poincare, Sur les equations lineaires aux differentielle
ordinaires et aux differences finies, 
{\it Amer. J. Math.} (1885) 203--258.
\newline \newline \vskip -.25in \noindent
[12] H. Sherkat, Investigation of the Hecke group $G_5$
and its Eisenstein series, (undergraduate thesis)
 http://halfaya.org/ucla/research/sherkat.pdf,
 U.C.L.A, Los Angeles, 2007.
\end{document}